\documentclass{article}

\usepackage{amsmath,amssymb}
\usepackage{amsthm}
\usepackage{mathrsfs}
\usepackage{xypic}
\xyoption{all}

\theoremstyle{plain}
\newtheorem{theorem}{Theorem}[section]
\newtheorem{lemma}[theorem]{Lemma}
\newtheorem{corollary}[theorem]{Corollary}
\newtheorem{proposition}[theorem]{Proposition}
\theoremstyle{definition}

\theoremstyle{remark}
\newtheorem{remark}[theorem]{Remark}

\DeclareMathOperator{\id}{id}

\DeclareMathOperator{\Spec}{Spec}
\DeclareMathOperator{\Hom}{Hom}
\DeclareMathOperator{\HOMs}{\mathscr{H}\!\!\mathscr{O}\!\!\mathscr{M}}
\DeclareMathOperator{\HOM}{HOM}
\DeclareMathOperator{\Pic}{\mathscr{P}\!\mathit{ic}}
\DeclareMathOperator{\Ext}{Ext}
\DeclareMathOperator{\Aut}{Aut}
\DeclareMathOperator{\Isom}{Isom}
\newcommand{\cosq}[2]{\mathrm{cosq}_0(#1 \! \to \! #2)}
\newcommand{\cosqzero}[1]{\mathrm{cosq}_0(#1)}

\DeclareMathOperator{\Ob}{Ob}

\newcommand{\defm}{\mathrm{Defm}}

\newcommand{\X}{\mathscr{X}} 
\newcommand{\Y}{\mathscr{Y}} 
\newcommand{\Z}{\mathscr{Z}} 
\newcommand{\T}{\mathscr{T}} 
\newcommand{\G}{\mathscr{G}} 
\newcommand{\F}{\mathscr{F}} 
\newcommand{\I}{\mathscr{I}} 
\renewcommand{\H}{\mathscr{H}} 

\renewcommand{\phi}{\varphi}
\renewcommand{\dot}{\bullet}

\newcommand{\ti}[1]{\widetilde{#1}}
\newcommand{\isomto}{\overset{\sim}{\longrightarrow}}
\newcommand{\tensor}{\otimes}
\newcommand{\ProjLim}{\mathop{\underleftarrow{\rm Lim}}}
\newcommand{\To}{\Rightarrow}
\newcommand{\Gm}{\mathbb{G}_m}

\title{Hom stacks}

\author{Masao Aoki\\
\normalsize
Department of Mathematics,
Kyoto University\\
\normalsize
\texttt{aoki@math.kyoto-u.ac.jp}}
\date{}

\begin{document}

\maketitle

\begin{abstract}
 We study Hom 2-functors
 parameterizing 1-morphisms of algebraic stacks,
 and prove that they are representable by algebraic stacks
 under certain conditions, using Artin's criterion. As an application
 we study Picard 2-functors which parameterize line bundles on
 algebraic stacks.
\end{abstract}

\section{Introduction}
\label{sec:intro}

Let $S$ be an affine noetherian scheme over an excellent Dedekind domain.
Let $\X$ and $\Y$ be algebraic stacks over $S$.
The Hom 2-functor $\HOMs(\X,\Y)$ is a contravariant 2-functor from
the category of affine noetherian schemes over $S$
to the 2-category of groupoids given by
\[
 \HOMs(\X,\Y)(T) = \HOM_T(\X \times_S T, \Y \times_S T).
\]
The right hand side is the groupoid of 1-morphisms.

The purpose of this paper is to show the following theorem:
\begin{theorem}
 \label{thm:main}
 If $\X$ is proper and flat over $S$
 and $\Y$ is of finite presentation over $S$,
 the 2-functor $\H=\HOMs(\X,\Y)$ is
 an algebraic stack in Artin's sense \cite{Ar}.
\end{theorem}
Here ``in Artin's sense'' means that the diagonal $\H \to \H \times_S
\H$ is representable and locally of finite type.
\\

It is already known (see \cite[4.1]{Ol.hom}) that if $X$ is
a proper flat algebraic space and
$Y$ is a separated algebraic space of finite type, the functor $\HOMs(X,Y)$ is
representable by an algebraic space. Moreover if $X$ and $Y$ are
quasi-projective schemes, $\HOMs(X,Y)$ is also a quasi-projective scheme.
This is proved by the fact that the map
\begin{eqnarray*}
  \HOMs(X,Y) &\to& \mathrm{Hilb}(X \times Y)\\
 f &\mapsto& \mbox{graph of $f$}
\end{eqnarray*}
is representable by an open immersion.

Unfortunately, we cannot use this technique in the case of algebraic
stacks, because we do
not have ``Hilbert stacks'' for algebraic stacks yet.
The $\mathrm{Quot}$ functors of Olsson and Starr
(\cite{OS},\cite{Ol.proper})  do not work for our purpose.
The functor $\mathrm{Quot}_{O_{\X \times \Y}}$ parameterizes closed substacks
of $\X \times \Y$, but graphs of 1-morphisms are not closed substacks
in general, even if the stacks $\X$ and $\Y$ are separated.
For instance, the graph of $\id: \X \to \X$ is the diagonal
$\X \to \X \times \X$, which is not a closed immersion unless $\X$ is
representable by an algebraic space.

Olsson \cite{Ol.hom} studied this problem when $\X$ and $\Y$ are
Deligne-Mumford stacks. He investigated the map
\[
 \HOMs(\X,\Y) \to \HOMs(\underline{\X},\underline{\Y})
\]
mapping a morphism to that of its coarse moduli spaces.
Even this technique does not work for Artin stacks,
because they do not have coarse moduli spaces in general.
\\

We prove Theorem \ref{thm:main} by verifying Artin's condition \cite{Ar}
directly. The most essential part of the proof is the deformation theory of
morphisms of algebraic stacks, based on the author's previous work
\cite{Ao}.

As an application, we prove that the Picard 2-functor \cite[14.4.7]{LM} that
parameterizes line bundles on an algebraic stack is representable by an
algebraic stack in Artin's sense.
This is a generalization of Artin's results on
algebraic spaces (\cite[7.3]{Ar1}, \cite[Appendix 2]{Ar}).

\subsection{Conventions and notations}

In this paper we refer to \cite{LM} for definitions and basic properties
of algebraic stacks. Especially we assume all algebraic stacks are
quasi-separated \cite[4.1]{LM} unless mentioned.
Algebraic stacks as in Artin's definition \cite[5.1]{Ar} are called
``algebraic stack in Artin's sense''.

We denote schemes and algebraic spaces by Italic
letters like $X,Y$ and $T$, and algebraic stacks
by script letters like $\X,\Y$ and $\T$.
Subscripts like $\X_T$ mean base change $\X \times_S T$.
Superscripts like $X^\dot$ are used to denote simplicial algebraic
spaces.

\subsection{Acknowledgements}
The author would like to express his thanks to Professor Fumiharu Kato
for valuable suggestions and advice on this paper,
and to Dr. Olsson, Mr. Iwanari and Dr. Yasuda
for useful comments and conversations.

Financial support is provided by Japan Society of Promotion of Science.

\section{Deformation of morphisms of algebraic stacks}

In this section we study the deformation theory of 1-morphisms
of algebraic stacks. This is a generalization of Illusie's work
\cite[III 2.2]{Il}.

\subsection{Definitions and Statements}

Deformations of 1-morphisms are defined as follows.
Let $\X$ and $\Y$ be separated algebraic stacks over a scheme $T$
and $f:\X \to \Y$ a 1-morphism over $T$.
Consider the 2-commutative diagram of solid arrows:
\[
 \xymatrix{
 \X \ar[rr]_{i} \ar[rd]^{f} \ar[rdd]& & \ti{\X} \ar[rdd] \ar@{.>}[rd]^{\ti{f}}\\
 & \Y \ar[rr]_{j} \ar[d] & & \ti{\Y} \ar[d]\\
 & T \ar[rr]_{k} & & \ti{T}.
}
\]
Here $i,j$ and $k$ are closed immersions defined by square-zero ideals
$I,J$ and $K$.
Then a deformation of $f$ is a pair $(\ti{f},\lambda)$ where
$\ti{f}$ is a 1-morphism from $\ti{\X}$ to $\ti{\Y}$ over $\ti{T}$
and $\lambda : \ti{f} \circ i \Rightarrow j \circ f$ is a 2-isomorphism.
A morphism from $(\ti{f}, \lambda)$ to $(\ti{g},\mu)$ is a 2-morphism $\alpha :
\ti{f} \Rightarrow \ti{g}$ such that the 2-morphisms
\[
 i^* \alpha \circ \mu, \lambda : \ti{f} \circ i \Rightarrow j \circ f
\]
are equal.

We denote the category of deformations of $f$ by $\defm_T(f)$
and the set of its isomorphism classes by $\overline{\defm}_T(f)$.

In this section we prove the following generalization of \cite[III 2.2.4]{Il}.
\begin{theorem}
 \label{thm:defm}
 $ $
 \begin{enumerate}
  \item There exists an obstruction $o \in \Ext^1(Lf^*L_{\Y/T},I)$
	whose vanishing is equivalent to the existence of a deformation.
  \item If $o=0$, the set $\overline{\defm}_T(f)$
	is a torsor under $\Ext^0(Lf^*L_{\Y/T},I)$.
  \item The automorphism group of any deformation of $f$ is isomorphic to\\
	$\Ext^{-1}(Lf^*L_{\Y/T},I)$.
 \end{enumerate}
\end{theorem}

In the proof of Theorem \ref{thm:defm}, we need the deformation theory
of morphisms of schemes over algebraic stacks.

Let $\T$ be an algebraic stack, $x:X \to \T$ and $y: Y \to \T$
schemes over $\T$, and $f:X \to Y$ a morphism of schemes with $y \circ f = x$.
Consider the diagram of solid arrows:
\[
 \xymatrix{
 X \ar[rr]_{i} \ar[rd]^{f} \ar[rdd]_{x}& &
 \ti{X} \ar[rdd]_(.3){\ti{x}} \ar@{.}[rd]^{\ti{f}}\\
 & Y \ar[rr]_{j} \ar[d]^{y} & & \ti{Y} \ar[d]^{\ti{y}}\\
 & \T \ar[rr]_{k} & & \ti{\T}.
}
\]
Here $i,j$ and $k$ are closed immersions defined by square-zero ideals
$I,J$ and $K$.
Then we define a deformation of $f$ to be a pair $(\ti{f},\gamma)$ where
$\ti{f}$ is a morphism $\ti{X} \to \ti{Y}$ which satisfies
$\ti{f} \circ i = j \circ f$ and $\gamma$ is a 2-isomorphism
$\ti{y} \circ \ti{f} \To \ti{x}$ whose restriction  $y \circ f \To x$
is equal to the identity.

We denote the set of deformations of $f$ by $\defm_\T(f)$.

\begin{proposition}
 \label{prop:defm_sch/ch}
 $ $
 \begin{enumerate}
  \item There exists an obstruction $o \in \Ext^1(Lf^*L_{Y/\T},I)$
	whose vanishing is equivalent to the existence of a deformation.
  \item If $o=0$, $\defm_T(f)$
	is a torsor under $\Ext^0(Lf^*L_{Y/\T},I)$.
 \end{enumerate}
\end{proposition}

\begin{remark}
 The torsor actions and isomorphisms in Theorem \ref{thm:defm} and
Proposition \ref{prop:defm_sch/ch} are functorial on $\X, \Y$ and $T$
 etc.
For example, if $T \to U$ is a morphism of schemes,
we have the natural ``forgetting'' map
\[
 C: \overline{\defm_T(f)} \to \overline{\defm_U(f)}
\]
and the group
homomorphism
\[
 D:\Ext^0(Lf^*L_{\Y/T},I) \to \Ext^0(Lf^*L_{\Y/U},I)
\]
induced by the morphism $L_{\Y/U} \to L_{\Y/T}$ \cite[17.3(3)]{LM}.
Then for any $[\ti{f}] \in \overline{\defm_T(f)}$ and $\sigma \in
\Ext^0(Lf^*L_{\Y/T},I)$, we have
\[
 C(\sigma \cdot [\ti{f}]) = D(\sigma) \cdot C([\ti{f}]).
\]
Note that this is true for schemes and simplicial algebraic spaces
(see the proof of \cite[III 2.2.4]{Il}).
We prove a special case of this for Proposition \ref{prop:defm_sch/ch}
which is necessary for the proof of Theorem \ref{thm:defm}.
A proof for the general case is straightforward.
\end{remark}

\subsection{Proof of Proposition \ref{prop:defm_sch/ch}}

The strategies of proofs of Theorem \ref{thm:defm} and Proposition
\ref{prop:defm_sch/ch} are the same as those of \cite{Ao} and \cite{Ol}.

\begin{description}
 \item[Step 1:] Choose good presentations of algebraic stacks
 and make associated
 simplicial algebraic spaces.
 \item[Step 2:] Compare deformations
 in the 2-category of algebraic stacks and
 those in the category of simplicial algebraic spaces.
 \item[Step 3:] Compare the $\Ext$ groups.
\end{description}

\begin{proof}[Proof of Proposition \ref{prop:defm_sch/ch}]
Let $\ti{P}^0 : \ti{T}^0 \to \ti{\T}$ be a presentation of $\ti{\T}$
and $T^0 = \ti{T}^0 \times_{\ti{\T}} \T$.
Then $P^0 : T^0 \to \T$ is a presentation of $\T$.
 Let $T^\dot = \cosq{T^0}{\T}$ and $\ti{T}^\dot
 = \cosq{\ti{T}^0}{\ti{\T}}$.
 Consider the diagram obtained by base changes $T^\dot \to \T$ and
 $\ti{T}^\dot \to \ti{\T}$:
 \[
 \xymatrix{
 X^\dot \ar[rr]_{i^\dot} \ar[rd]^{f^\dot} \ar[rdd]_{x^\dot}& &
 \ti{X}^\dot \ar[rdd]_(.3){\ti{x}^\dot} \ar@{.}[rd]^{\ti{f}^\dot}\\
 & Y^\dot \ar[rr]_{j^\dot} \ar[d]^{y^\dot} & &
 \ti{Y}^\dot \ar[d]^{\ti{y}^\dot}\\
 & T^\dot \ar[rr]_{k^\dot} & & \ti{T}^\dot.
 }
 \]
 Then by construction $\ti{X^\dot} \cong \cosq{\ti{X^0}}{\ti{X}}$ and
 $\ti{Y^\dot} \cong \cosq{\ti{Y^0}}{\ti{Y}}$.
 Therefore $\ti{f^\dot} : \ti{X^\dot} \to \ti{Y^\dot}$ descends to a morphism
 $\ti{f} : \ti{X} \to \ti{Y}$.
 Thus we can define a map $A':
 \defm_{T^\dot}(f^\dot) \to \defm_{\T}(f)$.
 
 The map $A'$ is bijective:
 the inverse is obtained by the base change.

 Let $I^\dot = \ker(O_{\ti{X}^\dot} \to O_{X^\dot})$.
 By the construction of the cotangent
 complex \cite[17.5]{LM}, the homomorphisms
 \[
 {P_X^\dot}^* : \Ext^i(Lf^*L_{Y/\T},I) \to
 \Ext^i({f^\dot}^*L_{Y^\dot/T^\dot}, I^\dot)
 \]
 are isomorphisms for all $i$.

 By \cite[III 2.2.4]{Il}, the obstruction for the existence of deformation of
 $f^\dot$ is in $\Ext^1({f^\dot}^*L_{Y^\dot/T^\dot}, I^\dot)$ and
 the set $\defm(f^\dot)$ is a torsor under
 $\Ext^0({f^\dot}^*L_{Y^\dot/T^\dot}, I^\dot)$. This proves the 
 proposition.
\end{proof}

Next we prove that the action of $\Ext$ groups are functorial on $\T$.

Let $f:X \to Y$ be a morphism over $\T$ as
in Proposition \ref{prop:defm_sch/ch} and
$\T \to U$ a morphism to a scheme.
Here we consider a deformation diagram:
\[
 \xymatrix@R1.2pc{
 X \ar[rd]^{f} \ar[rdd]^(.6){x} \ar@/_/[rddd] \ar[rr]& &
 \ti{X} \ar@{.>}[rd]^{\ti{f}} \ar[rdd]^(.7){\ti{x}} \ar@/_/[rddd]\\
 & Y \ar[d]^{y} \ar[rr] & & \ti{Y} \ar[d]^{\ti{y}}\\
 & \T \ar[d] \ar[rr] & & \ti{\T} \ar[d] \\
 & U \ar[rr] & & \ti{U}
}
\]

\begin{proposition}
\label{prop:defm_sch/ch.com}
The natural map
\[
 C:\defm_{\T}(f) \to \defm_{U}(f)
\]
is compatible with the homomorphism of groups
\[
 D: \Ext^0(Lf^*L_{Y/\T},I) \to \Ext^0(f^*L_{Y/U},I).
\]
\end{proposition}

\begin{proof}
 Let $T^\dot = \cosq{T^0}{\T}$ be the simplicial algebraic space as in
 the proof of Proposition \ref{prop:defm_sch/ch}.
 Consider the diagram obtained by base change:
\[
 \xymatrix@R1.2pc{
  & X^\dot \ar[ld]^{P_X^\dot} \ar[rrd]^{f^\dot} \ar[rrdd]\\
 X \ar[rrd]^{f} \ar[rrdd] & & & Y^\dot \ar[ld]^(.7){P_Y^\dot} \ar[d]\\
 & & Y \ar[d] & T^\dot \ar[ld]\\
 & & \T \ar[d]\\
 & & U.
 }
\]

The map $C$ factors as
\begin{gather*}
 \defm_{\T}(f) \stackrel{C_1}{\longrightarrow}
 \defm_{T^\dot}(f^\dot) \stackrel{C_2}{\longrightarrow}
 \defm_{U}(f^\dot)\\
 \stackrel{C_3}{\longrightarrow}
 \defm_U(P_Y^\dot \circ f^\dot) = \defm_U(f \circ P_X^\dot) 
 \stackrel{C_4}{\longrightarrow}
 \defm_U(f)
\end{gather*}
and $D$ factors as
\begin{gather*}
 \Ext^0(Lf^*L_{Y/\T},I) \stackrel{D_1}{\longrightarrow}
 \Ext^0({f^\dot}^*L_{Y^\dot/T^\dot}, I^\dot)
  \stackrel{D_2}{\longrightarrow}
 \Ext^0({f^\dot}^*L_{Y^\dot/U}, I^\dot)\\
  \stackrel{D_3}{\longrightarrow}
 \Ext^0((P_Y^\dot \circ f^\dot)^*L_{Y/U}, I^\dot)
 = \Ext^0((f \circ P_X^\dot)^*L_{Y/U}, I^\dot)\\
  \stackrel{D_4}{\longrightarrow}
 \Ext^0(f^*L_{Y/U},I).
\end{gather*}
 The compatibility of isomorphisms $C_1$ and $D_1$ is
 obvious by the definition of
 the action of $\Ext^0(Lf^*L_{Y/\T},I)$ in the proof of
 Proposition \ref{prop:defm_sch/ch}. That of $C_2$ and $D_2$ follows
 from the case of simplicial algebraic spaces.
 For $C_3$ and $D_3$, it follows from the definition of the morphism
 ${P_Y^\dot}^* L_{Y/U} \to L_{Y^\dot/U}$ \cite[II 1.2.7]{Il}.
 For $C_4$ and $D_4$, it is trivial.
\end{proof}

\subsection{Proof of Theorem \ref{thm:defm}: Step 1}

Let $P_Y: Y^0 \to \Y$ be a presentation of $\Y$, $\X' = \X \times_{\Y} Y^0$
and $X^0 \to \X'$ a presentation of $\X'$.
Then the composition $P_X : X^0 \to \X' \to \X$ is a presentation of
$\X$.
We may assume $X^0$ and $Y^0$ are
affine. Since $X^0 \to \X$ and $Y^0 \to \Y$ are smooth, we have the
unique deformations $\ti{X^0} \to \ti{\X}$ and $\ti{Y^0} \to \ti{\Y}$.
Let $X^\dot = \cosq{X^0}{\X}$ etc. We obtain the following diagram:
\[
 \xymatrix@R1.2pc{
 X^\dot \ar[dd]^{P_X^\dot} \ar[rd]^{f^\dot} \ar[rr] & &
 \ti{X^\dot} \ar[dd]_(.3){\ti{P_X^\dot}} \ar@{.}[rd]^{\ti{f^\dot}}\\
 & Y^\dot \ar[dd]^(.3){P_Y^\dot} \ar[rr] &
 & \ti{Y^\dot} \ar[dd]^{\ti{P_Y^\dot}} \\
 \X \ar[rd]^{f} \ar[rdd]_{x} \ar[rr]  & &
 \ti{\X} \ar@{.}[rd]^{\ti{f}} \ar[rdd]_(.3){\ti{x}}\\
 & \Y \ar[d]^{y} \ar[rr] &
 & \ti{\Y} \ar[d]^{\ti{y}}\\
 & T \ar[rr] &
 & \ti{T}.
}
\]
Let $I^\dot = \ker(O_{\ti{X}^\dot}
\to O_{X^\dot}) \cong {P_X^\dot}^* I$.

\subsection{Proof of Theorem \ref{thm:defm}: Step 2}

The map
\[
 A:\defm_T(f^\dot) \to \overline{\defm_T(f)}
\]
is defined by
sending $\ti{f}^\dot : \ti{X}^\dot \to \ti{Y}^\dot $ to the morphism
of associated stacks $\ti{f}: \ti{\X} \to \ti{\Y}$.

\begin{proposition}
\label{prop:A.surj}
The map $A$ is surjective. 
\end{proposition}

\begin{proof}
 Fix $[\ti{f}] \in \overline{\defm_T(f)}$. First we claim that
 $[\ti{f}]$ is in the image of $A$ if $\defm_{\Y}(f^0)$ is not empty.
 To see this, let $(\ti{f^0}, \gamma) \in \defm_{\Y}(f^0)$.
 We define $\ti{f}^\dot = \cosqzero{\ti{f}^0,\gamma}
 : \ti{X}^\dot \to \ti{Y}^\dot$
 as follows. Since $\ti{X}^\dot$ and $\ti{\Y}^\dot$ are the images of
 $\mathrm{cosq}$, by the similar discussion as in \cite[3.1.3]{Ao},
 to give $\ti{f^\dot}$ it suffices to
 give $\ti{f}^1:\ti{X^1} \to \ti{Y^1}$. This
 is equivalent to giving a triple $(\ti{f}^0 \circ p_1, \ti{f}^0 \circ
 p_2,\epsilon)$, where
\[
 \epsilon : P_Y \circ \ti{f}^0 \circ p_1
 \Rightarrow P_Y \circ \ti{f}^0 \circ p_2
\]
 is a 2-morphism.
 Now we put $\epsilon = p_2^* \gamma \circ p_1^* \gamma^{-1}$.
 Then $A(\ti{f}^\dot) = [\ti{f}]$.

 By Proposition \ref{prop:defm_sch/ch}
the obstruction for the existence of $(\ti{f^0},\gamma)$
is in\\ $\Ext^1(L{f^0}^*L_{Y^0/\Y},I^0)$. This group is zero because
$X^0$ is affine and $L_{Y^0/\Y}$ is quasi-isomorphic to a locally free
 sheaf $\Omega_{Y^0/\Y}$. 
\end{proof}

\begin{corollary}
 The obstruction for existence of a deformation of $f$ is in\\
 $\Ext^1({f^\dot}^* L_{Y^\dot/T}, I^\dot)$.
\end{corollary}

For each $[\ti{f}] \in \overline{\defm_T(f)}$, let $C$ be the
composition of maps
\[
 \defm_{\Y}(f^0) \stackrel{\mathrm{cosq}_0}{\isomto}
 \defm_{\Y}(f^\dot) \stackrel{\mbox{``forget''}}{\longrightarrow} \defm_T(f^\dot).
\] 
By Proposition \ref{prop:defm_sch/ch.com},
this is compatible with the group homomorphism
\[
 D: \Ext^0(L{f^0}^*L_{Y^0/\Y},I^0) \to \Ext^0({f^\dot}^*
 L_{Y^\dot/T},I^\dot).
\]

\begin{proposition}
 The set $\overline{\defm_T(f)}$ is the set of
 $\Ext^0(L{f^0}^*L_{Y^0/\Y},I^0)$-orbits in $\defm_T(f^\dot)$
 by the action induced by $D$.
\end{proposition}
\begin{proof}
 Suppose that $\ti{f}^\dot, \ti{g}^\dot \in \defm_T(f^\dot)$ satisfy
 $A(\ti{f^\dot}) = A(\ti{g^\dot}) = [\ti{f}]$. 
 Then there exists $(\ti{f^0},\gamma),(\ti{g^0},\delta)
 \in \defm_{\Y}(f^0)$ such that
 $C(\ti{f^0},\gamma) = \ti{f^\dot}$ and
 $C(\ti{g^0},\delta) = \ti{g^\dot}$.
 Since $\defm_{\Y}(f^0)$ is a
 $\Ext^0(L{f^0}^*L_{Y^0/\Y},I^0)$-torsor,
 there exists\newline
 $\sigma \in \Ext^0(L{f^0}^*L_{Y^0/\Y},I^0)$
 such that $\sigma \cdot (\ti{f^0},\gamma) = (\ti{g^0},\delta)$. Hence
 $D(\sigma) \cdot \ti{f^\dot} = \ti{g^\dot}$.

 Conversely, suppose that $\ti{f}^\dot, \ti{g}^\dot \in \defm_T(f^\dot)$
 satisfy $D(\sigma) \cdot \ti{f^\dot} = \ti{g^\dot}$ for some
 $\sigma \in \Ext^0(L{f^0}^*L_{Y^0/\Y},I^0)$. Let
 $[\ti{f}] = A(\ti{f^\dot})$ and choose $(\ti{f^0},\gamma) \in \defm_{\Y}(f^0)$
 such that $C(\ti{f^0},\gamma) = \ti{f^\dot}$. Then
 $C(\sigma \cdot (\ti{f^0},\gamma)) = D(\sigma) \cdot \ti{f^\dot}
 = \ti{g^\dot}$.
 Therefore $A(\ti{g^\dot}) = [\ti{f}]$.
\end{proof}

\begin{proposition}
 Fix an object $\ti{f}$ of $\defm_T(f)$.
 Then $\Aut(\ti{f})$, the group of automorphisms of deformations,
 is isomorphic to $\ker(D)$.
\end{proposition}

\begin{proof}
 Fix $\ti{f^\dot} \in \defm_T(f^\dot)$ such that $A(\ti{f^\dot}) = [\ti{f}]$
 and $(\ti{f}^0, \gamma) \in C^{-1}(\ti{f}^\dot)$.

 First we identify $\Aut(\ti{f})$ with a subset of $\defm_\Y(f^0)$
and construct set-theoretical bijection from $\Aut(\ti{f})$ to
$C^{-1}(\ti{f}^\dot)$.
Let $\alpha \in \Aut(\ti{f})$ and let $\beta$ be the composition of 2-morphisms
\[
 \ti{P_Y} \circ \ti{f^0}
 \stackrel{\gamma^{-1}}{\Longrightarrow}
 \ti{f} \circ \ti{P_X}
 \stackrel{\ti{P_X}^* \alpha}{\Longrightarrow}
 \ti{f} \circ \ti{P_X}
 \stackrel{\gamma^{-1}}{\Longrightarrow}
 \ti{P_Y} \circ \ti{f^0}.
\]
Then the triple $(\ti{f}^0,
 \ti{f}^0,\beta)$ defines a morphism
\[
 d_\alpha : \ti{X}^0 \to \ti{Y^0} \times_{\ti{\Y}} \ti{Y^0}
 = \ti{Y}^1.
\]
This is an element of $\defm_{Y^0}(\Delta \circ f^0)$.
Here $Y^1$ is a scheme over $Y^0$ by $p_1:Y^1 \to Y^0$.
\[
 \xymatrix{
 X^0 \ar[rr] \ar[rd]^{\Delta \circ f^0} \ar[rdd]_{f^0} & &
 \ti{X^0} \ar[rd]^{d_\alpha} \ar[rdd]_(.3){\ti{f^0}}\\
 & Y^1 \ar[rr] \ar[d]_(.3){p_1} & & \ti{Y^1} \ar[d]^{p_1}\\
 & Y^0 \ar[rr] \ar@<-1ex>@/_/[u]_{\Delta} & & \ti{Y^0}
}
\]
The map
\begin{eqnarray*}
 p_1^* : \defm_{\Y}(f^0) = \defm_{\Y}(p_1 \circ \Delta \circ f^0)
  &\to& \defm_{Y^0}(\Delta \circ f^0)\\
 ({\ti{f^0}}',\gamma') &\mapsto&
(\ti{f}^0,{\ti{f^0}}', {\gamma'}^{-1} \circ \gamma)
\end{eqnarray*}
is a bijection and compatible with the isomorphism
\begin{gather*}
  p_1^* : \Ext^0(L{f^0}^* L_{Y^0/\Y},I^0)
 = \Ext^0(L(p_1 \circ \Delta \circ f^0)^* L_{Y^0/\Y},I^0)\\
 \isomto
 \Ext^0(L(\Delta \circ f^0)^* L_{Y^1/Y^0}, I^0)
\end{gather*}
induced by $p_1^* L_{Y^0/\Y} \cong L_{Y^1/Y^0}$.

Now $({\ti{f^0}}',\sigma')$ is in $C^{-1}(\ti{f}^\dot)$ if and only
if $\ti{f^0}' = \ti{f}^0$ and $p_2^*\gamma' \circ
p_1^*{\gamma'}^{-1} = p_2^* \gamma \circ p_1^* \gamma^{-1}$.
The latter is equivalent to
\[
 p_1^*({\gamma'}^{-1} \gamma) = p_2^*({\gamma'}^{-1} \gamma),
\]
which implies the existence of $\alpha \in \Aut(\ti{f})$
such that $\gamma' \circ \gamma^{-1} = \gamma \circ P_X^* \alpha \circ
\gamma^{-1}$.

Thus we can identify $\Aut(\ti{f})$ with $C^{-1}(\ti{f}^\dot)$ as subsets
of $\defm_{\Y}(f^0)$.
\\

 Next we see that the group structure of $\Aut(f)$
 is compatible with that of $\ker(D)$ acting on $C^{-1}(\ti{f}^\dot)$.
 The composition $\alpha \circ
 \alpha'$ corresponds to the morphism
 \[
  d_{\alpha \circ \alpha'} =
 (\ti{f^0},\ti{f^0},\gamma \circ \ti{P_X}^* \alpha \circ
 \ti{P_X}^* \alpha' \circ \gamma^{-1}) : \ti{X^0} \to \ti{Y^1}.
 \]
 This is equal to the composition
\[
 \ti{X}^0 \stackrel{(d_{\alpha'},d_{\alpha})}{\longrightarrow}
 \ti{Y}^1 \times_{p_1 \ti{Y}^0 p_2} \ti{Y}^1
 = \ti{Y}^2 \stackrel{p_{13}}{\longrightarrow} \ti{Y}^1.
\]
\[
 \xymatrix@C5pc{
 & \ti{Y^2} = \ti{Y^0} \times_{\ti{\Y}} \ti{Y^0} \times_{\ti{\Y}}
 \ti{Y^0} \ar[d]|{p_{13}} \ar@<-4ex>[d]|{p_{12}}
 \ar@<+4ex>[d]|{p_{23}}\\
 \ti{X^0} \ar[ur]^{(d_{\alpha'},d_{\alpha})}
 \ar[r]^{d_{\alpha}, d_{\alpha'}}_{d_{\alpha \circ \alpha'}}
 \ar[dr]_{\ti{f^0}} &
 \ti{Y^1} = \ti{Y^0} \times_{\ti{\Y}} \ti{Y^0}
 \ar@<+2ex>[u] \ar@<-2ex>[u] \ar@<-2ex>[d]|{p_1} \ar@<+2ex>[d]|{p_2}\\
 & \ti{Y^0} \ar[u]|{\Delta}
 }
\]
\\

On the other hand, the group structure of
\[
 \Ext^0((\Delta \circ f^0)^*L_{Y_1/Y_0},I^0)
 \cong \mathrm{Der}_{O_{Y^0}}(O_{Y^1},I^0)
\]
is given by taking sums of derivations $D_{\alpha},D_{\alpha'}:
O_{Y^1} \to I^0$ in the topos of \'{e}tale sheaves.
Pulling back by $p_{12}:Y^2 \to Y^1$,
we identify $D_{\alpha}$ with a derivation
\[
\begin{matrix}
 O_{Y^2} = O_{Y^1} \tensor_{p_1^* O_{Y^0} p_2^*} O_{Y^1}
 &\stackrel{D_{\alpha}}{\longrightarrow} & I^0\\
 x \tensor y &\mapsto & D_{\alpha}(x \tensor y) = x D_{\alpha}(1 \tensor y).
\end{matrix}
\]
in $\mathrm{Der}_{O_{Y^1}}(O_{Y^2},I^0)$.
Pulling back by $p_{23}:Y^2 \to Y^1$, $D_{\alpha'}$ is
identified with
\[
\begin{matrix}
 O_{Y^2} = O_{Y^1} \tensor_{p_1^* O_{Y^0} p_2^*} O_{Y^1}
 &\stackrel{D_{\alpha'}}{\longrightarrow}  & I^0\\
 x \tensor y &\mapsto & D_{\alpha'}(x \tensor y) = y D_{\alpha'}(x \tensor 1).
\end{matrix}
\]
The morphism $(d_{\alpha'},d_{\alpha})$ as above corresponds to
a derivation
\[
\begin{matrix}
 O_{Y^3} = O_{Y^1} \tensor_{p_1^* O_{Y^0} p_2^*} O_{Y^1} \tensor_{p_1^* O_{Y^0} p_2^*} O_{Y^1}
 &\stackrel{D}{\longrightarrow} &I^0\\
 x \tensor y \tensor 1 &\mapsto & y D_{\alpha'}(x \tensor 1)\\
 1 \tensor y \tensor z &\mapsto & y D_{\alpha}(1 \tensor z).
\end{matrix}
\]
Then the morphism $d_{\alpha \circ \alpha'}$ corresponds to the composition:
\[
\begin{matrix}
  O_{Y^2} = O_{Y^1} \tensor_{p_1^* O_{Y^0} p_2^*} O_{Y^1}
 & \stackrel{p_{13}^*}{\to} & O_{Y^3}
 & \stackrel{D}{\longrightarrow} & I^0\\
 x \tensor y  &\mapsto & x \tensor 1 \tensor y
 &\mapsto & D((x \tensor 1 \tensor 1)(1 \tensor 1 \tensor y))\\
 & & & & = yD(x \tensor 1 \tensor 1) + xD(1 \tensor 1 \tensor y)\\
 & & & & = D_{\alpha}(x \tensor y) + D_{\alpha'}(x \tensor y)
\end{matrix}
\]
Thus group structures of $\Aut(\ti{f})$ and
$\mathrm{Der}_{O_{Y^0}}(O_{Y^1},I^0)$ are compatible.
\end{proof}
 
\subsection{Proof of Theorem \ref{thm:defm}: Step 3}

The following lemma completes the proof of Theorem \ref{thm:defm}.
\begin{lemma}
 $ $
 \begin{enumerate}
  \item There is an isomorphism
	\[
	\Ext^1({f^\dot}^*L_{Y^\dot/T},I^\dot) \isomto
	\Ext^1(Lf^*L_{\Y/T},I).
	\]
  \item The cokernel of $D:\Ext^0(L{f^0}^*L_{Y^0/\Y},I^0) \to
	\Ext^0({f^\dot}^* L_{Y^\dot/T},I^\dot)$ is isomorphic to
	$\Ext^0(Lf^*L_{\Y/T},I)$.
  \item The kernel of $D$ is isomorphic to $\Ext^{-1}(Lf^*L_{\Y/T},I)$.
 \end{enumerate}
\end{lemma}

\begin{proof}
 The morphisms
\[
 Y^\dot \to \Y \to T
\]
 induce a triangle in $D(O_{Y^\dot})$
 \[
  L{P_Y^\dot}^*L_{\Y/T} \to L_{Y^\dot/T} \to L_{Y^\dot/\Y} \to
 L{P_Y^\dot}^*L_{\Y/T}[1],
 \]
 and this in turn induces a long exact sequence
\[
 \begin{matrix}
 & & & & 0 &\to &\Ext^{-1}(L{f^\dot}^*L{P_Y^\dot}^*L_{\Y/T},I^\dot)\\
 &\to &\Ext^{0}(L{f^\dot}^*L_{Y^\dot/\Y},I^\dot)
 &\to &\Ext^{0}({f^\dot}^*L_{Y^\dot/T},I^\dot)
 &\to &\Ext^{0}(L{f^\dot}^*L{P_Y^\dot}^*L_{\Y/T},I^\dot)\\
 &\to &\Ext^{1}(L{f^\dot}^*L_{Y^\dot/\Y},I^\dot)
 &\to &\Ext^{1}({f^\dot}^*L_{Y^\dot/T},I^\dot)
 &\to &\Ext^{1}(L{f^\dot}^*L{P_Y^\dot}^*L_{\Y/T},I^\dot)\\
 &\to &\Ext^{2}(L{f^\dot}^*L_{Y^\dot/\Y},I^\dot)
 &\to & \cdots
\end{matrix}
\]
By the similar discussion as in \cite[4.7]{Ol},
\[
 \Ext^i(L{f^\dot}^*L_{Y^\dot/\Y},I^\dot) \cong
 \Ext^i(L{f^0}^*L_{Y^0/\Y},I^0)
\]
and the right hand side is zero for $i > 0$.
The isomorphism ${P_X^\dot}^* : D^{+}(O_{\X}) \to D^{+}(O_{X^\dot})$
induces isomorphisms
\[
 \Ext^{i}(L{f^\dot}^*L{P_Y^\dot}^*L_{\Y/T},I^\dot)
 \cong \Ext^{i}(L{P_X^\dot}^*Lf^*L_{\Y/T},I^\dot)
 \cong \Ext^{i}(Lf^*L_{\Y/T},I).
\]
\end{proof}

\section{Artin's criterion}

In this section we prove Theorem \ref{thm:main} by
verifying the following Artin's criterion \cite[5.3]{Ar}.

\begin{enumerate}
 \item $\H$ is a limit-preserving stack.
 \item $\H$ satisfies Schlessinger's conditions.
       \begin{itemize}
	\item[(S1)] If $A' \to A$ and $B \to A$ are homomorphisms of
		    noetherian rings over $S$
		    and $A' \to A$ is a small extension, then for any $f \in
		    \H(A)$ the natural functor
		    \[
		    \H_f(A' \times_A B) \to \H_f(A') \times \H_f(B)
		    \]
		    is an equivalence of categories.
		    Here $\H_f(R)$ denotes the subcategory of $\H(R)$
		    consisting of objects $g$ such that $g|_A \simeq f$
		    and morphisms $\alpha$ such that $\alpha|_A
		    = \id_f$.
	\item[(S2)] If $M$ is a finite $A$-module and $f \in \H(A)$,
		    then
		    \[
		     D_f(M) = \Ob \H_f(A+M) / \sim
		    \]
		    is a finite
		    $A$-module.
       \end{itemize}
 \item Compatibility with completion.\\
       If $A$ is a complete local noetherian ring with maximal ideal
       $m$, the functor
       \[
	\H(A) \to \ProjLim_n \H(A/m^{n+1})
       \]
       is an equivalence.
 \item Conditions on modules of obstruction, deformations and
       infinitesimal automorphisms.\\
       For any $f \in \H(A)$ and a finite $A$-module $M$,
       there exists a module of obstructions
       $O_f(M)$, a module of deformations $D_f(M)$ and a module of
       infinitesimal automorphisms $\Aut_f(M)$ which satisfy the
       following conditions:
       \begin{enumerate}
	\item compatibility with \'{e}tale localization:\\
	      If $A \to B$ is \'{e}tale and $g$ is a image of
	      $f$ in $\H(B)$,
	      \[
	      D_g(M \tensor B) \cong D_f(M) \tensor_A B
	      \]
	      etc.
	\item compatibility with completion:\\
	      If $m$ is a maximal ideal of $A$ and $\hat{A}$ is
	      a completion with respect to $m$,
	      \[
	       D_f(M) \tensor \hat{A} \cong 
	      \varprojlim D_f(M/m^n M)
	      \]
	      etc.
	\item constructibility:\\
	      There is a open dense set of points of finite type
	      $A \to k(p)$ such that
	      \[
	       D_f(M) \tensor k(p) \cong D_f(M \tensor k(p)).
	      \]
	      etc.
       \end{enumerate}
 \item For any $f \in \H(A)$ and $\alpha \in \Aut(f)$, if $\alpha |_k =
       \id$ for dense set of points of finite type $A \to k$, then
       $\alpha = \id$.
\end{enumerate}

\subsection{Preliminaries}

We can reduce many properties of $\H$ to that of $\Y$
by the following observations.

\begin{lemma}
 \label{lem:stack_gr}
 Let $\X$ and $\Y$ be algebraic stacks over $S$ and
 $X \to \X$ an epimorphism (e.g. a presentation of $\X$).
 Let $X^1 = X^0 \times_{\X} X^0$. Then the category
 $\HOM_S(\X,\Y)$ is equivalent to the following category:
 \begin{itemize}
  \item An object is a pair $(f^0,\alpha)$ where
	$f^0$ is an object of $\Y(X^0)$ and $\alpha : p_1^* f^0 \To p_2^* f^0$
	is a morphism in $\Y(X^1)$. 
  \item A morphism from $(f^0,\alpha)$ to $(g^0,\beta)$ is a
	morphism $\gamma : f^0 \To g^0$ in $\Y(X^0)$ such that
	$p_2^* \gamma \circ \alpha = \beta \circ p_1^* \gamma$
	in $\Y(X^1)$.
 \end{itemize}
\end{lemma}
\begin{proof}
 This follows immediately from the fact that $\X$ is a stack
 associated to the groupoid $X^1 \rightrightarrows X^0$
 by \cite[3.8]{LM}.
\end{proof}

\begin{lemma} 
 \label{lem:base_ch}
 Let $y:\Y \to S$ be an algebraic stack over a scheme $S$,
 $\phi: T \to S$ a morphism
 of schemes and $x: \X_T \to T$ an algebraic stack over $T$. Then the natural
 functor 
\[
 \HOM_T(\X_T,\Y_T) \to \HOM_S(\X_T,\Y)
\]
is an equivalence of categories.
\end{lemma}

\begin{proof}
If $\X_T$ is a scheme, this is clear by the construction of fiber products
\cite[2.2.2]{LM}.
In the general case, let $X^0 \to \X_T$ be a presentation
and $X^1 = X^0 \times_{\X} X^0$.
Then by the case of schemes we have
\begin{align*}
 \Y_T(X^0) &\simeq \Y(X^0)\\
 \Y_T(X^1) &\simeq \Y(X^1).
\end{align*}
The result follows from Lemma \ref{lem:stack_gr}.
\end{proof}

\subsection{Limit preserving stack}

Fix a presentation $X^0 \to \X$ and let $X^1 = X^0 \times_{\X} X^0$.
Then if $\{U_i \to U\}$ is an \'{e}tale covering,
so is $\{X^k_{U_i} \to X^k_{U}\}$ for $k=0,1$.
The conditions of stacks for
$\H$ follows from those of $\Y$:
\begin{enumerate}
 \item Let $f$ and $g$ be objects of $\H(U)$ and
       $\phi, \psi : f \To g$ be morphisms in $\H(U)$.
       Suppose that $\phi|_i = \psi|_i$ in $\H(U_i)$ for all
       $i$.
       By Lemma \ref{lem:base_ch}, $\phi$ and $\psi$ are identified
       with morphisms in $\HOM(\X_U,\Y)$.
       Let $\phi'$ and $\psi'$, morphisms in $\Y(X^0_{U})$
       corresponding to $\phi$ and $\psi$ by Lemma \ref{lem:stack_gr}.
       Then $\phi'|_{X^0_{U_i}} = \psi'|_{X^0_{U_i}}$ for all $i$
       imply $\phi' = \psi'$. Hence $\phi = \psi$.

 \item Let $f$ and $g$ be objects of $\H(U)$ and
       $\phi_i : f|_i \To g|_i$ morphisms in $\H(U_i)$.
       Suppose that $\phi_i|_{ij} = \phi_j|_{ij}$ for all $i$
       and $j$. Let $(f^0,\alpha)$ and $(g^0,\beta)$ be pairs
       corresponding to $f$ and $g$, and
       $\phi'_i$ morphisms in $\Y(X^0_{U_i})$
       corresponding to $\phi_i$. Then $\phi'_i|_{X^0_{U_{ij}}}
       = \phi'_j|_{X^0_{U_{ij}}}$
       imply existence of $\psi' : f^0 \To g^0$ in $\Y(X^0_{U})$ such that
       $\psi'|_{X_{U_i}} = \phi'_i$.
       Since
       \[
       p_2^* \psi'|_{X_{U_i}} \circ \alpha|_{X_{U_i}}
       = \beta|_{X_{U_i}} \circ p_1^* \psi'|_{X_{U_i}}
       \]
       hold for all $i$,
       \[
	p_2^* \psi' \circ \alpha = \beta \circ p_1^* \psi'
       \]
       and $\psi'$ corresponds to a morphism
       $\psi : f \To g$ in $\H(U)$ such that $\psi|_i = \phi_i$.

 \item Let $f_i$ be objects of $\H(U_i)$ and $\phi_{ij}:
       f_i|_{ij} \To f_j|_{ij}$ morphisms in $\H(U_{ij})$ which
       satisfy cocycle conditions:
       \[
	\phi_{jk}|_{ijk} \circ \phi_{ij}|_{ijk}
       = \phi_{ik}|_{ijk}.
       \]
       Let $(f^0_i,\alpha_i)$ be pairs corresponding to
       $f_i$ and $\phi'_{ij}$ morphisms in $\Y(X^0_{U_{ij}})$ corresponding to
       $\phi_{ij}$.
       Then by the cocycle conditions
       \[
	\phi'_{jk}|_{X^0_{U_ijk}} \circ \phi'_{ij}|_{X^0_{U_ijk}}
       = \phi'_{ik}|_{X^0_{U_ijk}},
       \]
       there exists an object $f^0$ of $\Y(X^0_U)$ and
       morphisms $\psi'_i : f^0|_{X^0_{U_i}} \To f^0_i$
       such that $\phi'_{ij} \circ \psi'_i|_{X^0_{U_ij}}
       = \psi'_j|_{X^0_{U_ij}}$.
       Let
       \[
	\beta_i = p_2^*{\psi'_i}^{-1} \circ \alpha_i \circ p_1^*\psi'_i
       : p_1^*f^0|_{X^1_{U_i}} \To p_2^*f^0|_{X^1_{U_i}}.
       \]
       Then
       \begin{align*}
	\beta_i|_{X_{U_{ij}}} &=
	p_2^* {\psi'_i}^{-1}|_{X_{U_{ij}}} \circ \alpha_i|_{X_{U_{ij}}}
	\circ p_1^* \psi'_i|_{X_{U_{ij}}}\\
	&= p_2^* {\psi'_i}^{-1}|_{X_{U_{ij}}} \circ 
	p_2^* {\phi'_{ij}}^{-1} \circ
	\alpha_j|_{X_{U_{ij}}} \circ p_1^* \phi'_{ij}
	\circ p_1^* \psi'_i|_{X_{U_{ij}}}\\
	&= p_2^* {\psi'_j}^{-1}|_{X_{U_{ij}}} \circ \alpha_j|_{X_{U_{ij}}}
	\circ p_1^* \psi'_j|_{X_{U_{ij}}}\\
	&= \beta_j|_{X_{U_{ij}}}.
       \end{align*}
       Therefore there exists $\beta : p_1^* f^0 \To p_2^* f^0$ in
       $\Y(X^1_U)$ such that $\beta|_{X_{U_i}} = \beta_i$.
       The pair $(f^0,\beta)$ defines an object $f$ of $\H(U)$.
       The morphism $\psi'_i$ satisfies
       \[
       p_2^* \psi'_i \circ \beta|_{X_{U_i}} = \alpha_i \circ p_1^* \psi'_i.
       \]
       Therefore $\psi'_i$ corresponds to $\psi_i : f|_i \To f_i$
       such that $\phi_{ij} \circ \psi_i|_{ij} = \psi_j|_{ij}$.
  \end{enumerate}

$\H$ is limit-preserving by \cite[4.18]{LM}.

\subsection{Schlessinger's conditions}

First, let $\phi : A' \to A$ and $\psi : B \to A$ be
homomorphisms of noetherian rings over $S$
and suppose $\phi$ is a small extension.
Let $f \in \H(A)$.
By Lemma \ref{lem:base_ch}, the condition (S1') on $\H$ is equivalent to
the equivalence
\[
 \HOM_f(\X_{A' \times_A B}, \Y) \isomto
 \HOM_f(\X_{A'},\Y) \times \HOM_f (\X_B,\Y).
\]

Let $X^0 \to \X$ be a presentation. Since $\X$ is of finite type over
noetherian base, we may assume $X^0$ is a noetherian affine scheme $\Spec R$.

\begin{lemma}
 \label{lem:prod_tensor}
 The homomorphism
\begin{eqnarray*}
   \pi : R \tensor (A' \times_A B) &\to&
   (R \tensor A') \times_{R \tensor A} (R \tensor B)\\
 r \tensor (a', b) &\mapsto& (r \tensor a', r \tensor b)
\end{eqnarray*}
is an isomorphism.
\end{lemma}
\begin{proof}
 The kernel of the projection $A' \times_A B \to B$ is isomorphic to
 $\ker \phi$ and the kernel of $(R \tensor A') \times_{R \tensor A}
(R \tensor B) \to R \tensor B$ is isomorphic to $\ker(\id_R \tensor
 \phi)$. Since $R$ is flat, the horizontal sequences of the 
following diagram  are exact:
\[
\xymatrix@R1pc@C1pc{
 0 \ar[r] & R \tensor \ker \phi \ar[r] \ar@{=}[d] &
 R \tensor (A' \times_A B) \ar[r] \ar[d]_{\pi} &
 R \tensor B \ar[r] \ar@{=}[d]  & 0\\
 0 \ar[r] & R \tensor \ker \phi \ar[r] &
 (R \tensor A') \times_{R \tensor A} (R \tensor B) \ar[r] &
 R \tensor B \ar[r]& 0.
}
\]
It is easy to check that this diagram commutes. Therefore $\pi$ is 
an isomorphism.
\end{proof}

Let $X^1 = X^0 \times_\X X^0$ and
$(f^0,\alpha)$ a pair
correspond to $f:\X \to \Y$ as in Lemma \ref{lem:stack_gr}.
By the condition (S1') for $\Y$ and Lemma \ref{lem:prod_tensor},
we have an equivalence
\[
 \Y_{f^0}(X^0_{A' \times_A B}) \isomto
  \Y_{f^0}(X^0_{A'}) \times \Y_{f^0}(X^0_B)
\]
Since the functor $\Isom(p_1^*f^0, p_2^* f^0)$ is represented by an
algebraic space, we also have
\[
 \Isom_\alpha(p_1^*f^0_{X^0_{A' \times_A B}},p_2^* f^0_{X^0_{A' \times_A B}})
 \isomto
 \Isom_\alpha(p_1^*f^0_{X^0_{A'}}, p_2^* f^0_{X^0_{A'}})
 \times
 \Isom_\alpha(p_1^*f^0_{X^0_{B}}, p_2^* f^0_{X^0_{B}})
\]
These equivalences proves (S1').

By Theorem \ref{thm:main}, we have
\[
 D_{f_{X_0}}(M) \cong \Ext^0(Lf_{A_0}^* L_{\X_{A_0}/A_0}, x_{A_0}^* M).
\]
This is a finite $A_0$ module because $Lf_{A_0}^* L_{\X_{A_0}/A_0}$ is
coherent and $\X_{A_0}$ is proper over $A_0$.
This proves (S2).

\subsection{Compatibility with completion}
\label{ss:completion}

Let $A_n = A/m^{n+1}$. The functor
\[
 \H(A) \to \ProjLim \H(A_n).
\]
is equal to the functor
\[
 \pi : \HOM_A(\X_A,\Y_A) \to \ProjLim \HOM_{A_n}(\X_{A_n},\Y_{A_n}).
\]
Note that $\pi$ is a bijection if $\X$ and $\Y$  are schemes
\cite[5.4.1]{EGA}.
\\

First we reduce the problem to the case $\X_A$
is representable by a scheme $X_A$.

Consider the functor
\[
 \HOM(\X_A,\Y_A) \to \HOM((\X_A)_{\mathrm{red}},\Y_A).
\]
By Theorem \ref{thm:defm}, fibers of this functors are
described by $\Ext$ groups of the cotangent complexes.
They are isomorphic to the limits of those of the reductions
by the Grothendieck existence theorem for Artin stacks
\cite[8.1]{Ol.sheaves}. So we may suppose that $\X_A$ is
reduced.

Let $X^0_A \to \X_A$ be a proper surjection from a scheme
\cite[1.1]{Ol.proper}.
Since $\X_A$ is reduced, the surjective morphism $X^0_A \to \X_A$ is an
epimorphism.
By Lemma \ref{lem:stack_gr}, the functor $\pi$ is an equivalence if
the categories $\HOM(X^0_A,\Y_A)$ and $\HOM(X^1_A,\Y_A)$ are
equivalent to the limits of the reductions.\\

To see $\pi$ is fully faithful, let $f$ and $g$ be objects of the left hand
side. The functor $I=\Isom(f,g)$ is representable by a separated
algebraic space of finite type over $A$. So it suffices to show
the map
\[
 \pi' : \Hom(X_A,I) \to \varprojlim \Hom(X_{A_n},I)
\]
is bijective.

This map is surjective by the same argument as in
\cite[III 5.4.1]{EGA} using the Grothendieck existence theorem for algebraic
spaces \cite[V 6.3]{Kn}.

To see the injectivity of $\pi'$, let $\alpha$ and $\beta$ be
the elements of $\Hom(X_A,I)$.
The functor $I' = \Isom(\alpha,\beta)$ is representable
by a closed subscheme of $X_A$, and $\pi'$ is injective if
the map
\[
 \Hom(X_A,I') \to \varprojlim \Hom(X_{A_n},I')
\]
is surjective. This follows from \cite[III 5.4.1]{EGA}.
\\

To see $\pi$ is essentially surjective,
let $\{f_n\}$ be an object of the right hand side.

For each $n$, let $\G_n$ be the essential image \cite[3.7]{LM} of the morphism
$(\id,f_n): X_{A_n} \to X_{A_n} \times \Y_{A_n}$. More precisely,
for any scheme $T$ over $A_n$, the set of objects of $\G_n(T)$
is equal to $X_{A_n}(T)$ and the automorphisms group of an object $x$
is equal to the automorphism group of $f(x)$ in $\Y_{A_n}(T)$.

Then $\G_n$ is a closed substack of $X_{A_n} \times Y_{A_n}$,
and proper over $A_n$ since the composition $\G_n
\hookrightarrow X_{A_n} \times \Y_{A_n} \to X_{A_n}$ is proper.
Hence it corresponds to an ideal sheaf $\I_n$ whose support is proper
over $A_n$. By the Grothendieck
existence theorem for Artin stacks \cite[1.5]{Ol.proper},
there exists an ideal sheaf $\I$ of $\X_A$ 
with proper support whose reduction
on $\X_{A_n}$ is isomorphic to $\I_n$. Let $\G$ be the closed substack
of $X_A \times \Y_A$ corresponding to $\I$.
The stack $\G$ is proper over $A$.
Let $p : \G \to X_A$ the
composition $\G \hookrightarrow X_A \times \Y_A \to X_A$.

We claim that $p$ is an epimorphism. This follows from the following lemma.

\begin{lemma}
 \label{lem:reduction_epi}
\begin{enumerate}
 \item Let $Z$ and $T$ be proper algebraic spaces over $A$
       and $g:Z \to T$ a morphism over $A$.
       If all reductions of $g$ are isomorphisms (resp.\ closed immersions),
       then $g$ is an isomorphism (resp.\ a closed immersion).
 \item Let $\Z$ and $\T$ be proper algebraic stacks over $A$
       and $g:\Z \to \T$ a morphism over $A$.
       If all reductions of $g$ are epimorphisms, then $g$ is an 
       epimorphism.
\end{enumerate}
\end{lemma}
\begin{proof}
 As in \cite[I 4.6.8]{EGA}, we may suppose $\T=\Spec A$.
\begin{enumerate}
 \item The open subscheme of scheme-like points \cite[II 6.6]{Kn}
       contains the closed subscheme $\Z_{A_0}$. Therefore $U$ is
       equal to $\Z$ and $\Z$ is an scheme. The desired results follow from
       \cite[I 4.6.8]{EGA}. 
 \item By the decomposition of $g$ into an epimorphism and a monomorphism
       \cite[3.7]{LM}, it suffices to show that if all 
       reductions of $g$ are isomorphisms, so is $g$.
       
       Now it suffices to show that $\Z$ is an algebraic space.
       Consider the diagonal map
       \[
       \Delta : \Z \to \Z \times \Z.
       \]
       This is proper, representable and all its reductions are
       closed immersions. Therefore $\Delta$ is a closed immersion,
       which means $\Z$ is an algebraic space.
\end{enumerate}
\end{proof}

Now the category of morphisms from $X$ to $\G$ is equivalent to the
category of morphisms from the groupoid $\G \times_X \G
\rightrightarrows \G$ to $\G$. Construct a morphism $F : \G \to \G$
as follows. For any scheme $U$ and an object $x$ of $\G(U)$,
$F(x) = x$, and for any automorphism $\sigma$ of $x$,
$F(\sigma) = \id_x$.

For each $n$, the reduction $F_n:\G_n \to \G_n$ of $F$ factors through
$X_{A_n}$, hence gives a 2-isomorphism $\alpha_n :
p_1^* F_n \to p_2^* F_n$ in $\HOM(\G_n \times_{X_{A_n}} \G_n, \G_n)$.
Since the reduction is fully faithful, there exists
$\alpha : p_1^* F \to p_2^* F$ in $\HOM(\G \times_{X_A} \G, \G)$.
Therefore $F$ factors through $X_A$.

The composition
\[
 X_A \to \G \hookrightarrow X_A \times \Y_A \to \Y_A
\]
is the desired morphism.

\begin{remark}
 This discussion will be clearer if we use the theory of ``formal
 algebraic stacks'' \cite{Iw} by Iwanari.
\end{remark}

\subsection{Conditions on modules}

By Theorem \ref{thm:defm}, the modules $O_f(M)$, $D_f(M)$
and $\Aut_f(M)$ are represented as follows:
\begin{align*}
 O_f(M) &= \Ext^1(Lf^*L_{Y_A/A},x_A^* M)\\
 D_f(M) &= \Ext^0(Lf^*L_{Y_A/A},x_A^* M)\\
 \Aut_f(M) &= \Ext^{-1}(Lf^*L_{Y_A/A},x_A^* M)
\end{align*}
Here $x_A$ denotes the structural morphism $\X_A \to \Spec A$.

The compatibility with \'{e}tale localization is equivalent to
that the maps
\[
 \Ext^i(Lf^*L_{\Y_B/B},I \tensor B) \to
 \Ext^i(Lf^*L_{\Y_A/A},I) \tensor B \quad
(i=-1,0,1)
\]
are isomorphisms for any \'{e}tale localization $A \to B$.
Since $L_{B/A}=0$, we have $L_{\Y_B/B} \cong L_{\Y_A/A}$, which induces
the desired isomorphisms.

The compatibility with completion follows from \ref{ss:completion}.

The constructibility of these modules follows from
the semicontinuity theorem for proper algebraic stacks
(Theorem \ref{thm:semi-con}).

\subsection{Quasi-separation of the diagonal}

Let $f \in \H(A)$, $\alpha \in \Aut(f)$ and suppose that $\alpha |_k
= \id$ for a dense set of points $A \to k$. Fix a presentation $P : X^0
= \Spec R \to \X$. Then $P^* \alpha$ is an automorphism of
${P_A}^* f \in \Y(X^0_A)$. The set of points $R \tensor A \to k'$
which factors through $R \tensor k$ with $\alpha |_k = \id$ is dense
in $X^0_A$, and $P^* \alpha|_{k'} = \id$ on such points.
Hence $P^* \alpha = \id$ because $\Y$ is a quasi-separated
stack. This implies $\alpha= \id$.

\section{A remark on quasi-separation}

It is hard to show that the stack $\H$ is quasi-separated,
in other words, it is an algebraic stack in the sence of \cite[4.1]{LM}.
In the case of Deligne-Mumford stacks, Olsson \cite{Ol.hom} needed
some extra hyposeses to prove this.
In our case we have the following partial result.

\begin{proposition}
 Let $\X$ and $\Y$ as in Theorem \ref{thm:main}.
 Suppose that $\X=X$ is representable by an algebraic space
 and $\Y$ has a proper presentation $Y^0 \to \Y$.
 Then the stack $\H = \HOMs(X,\Y)$ is quasi-separated.
\end{proposition}

\begin{proof}
 What we have to show is that 
 if $f$ and $g$ are objects of $\HOMs(X,\Y)(T)$, 
 then the algebraic space $\Isom_T(f,g)$ is separated and quasicompact
 over $T$.
 
 Let $X_f = X_T \times_{f \Y_T} Y^0_T$,
 $X_g = X_T \times_{g \Y_T} Y^0_T$, $X^0_T = X_f \times_{X_T} X_g$ and
 $f^0,g^0 : X^0_T \to Y^0_T$ morphisms induced by $f$ and $g$.
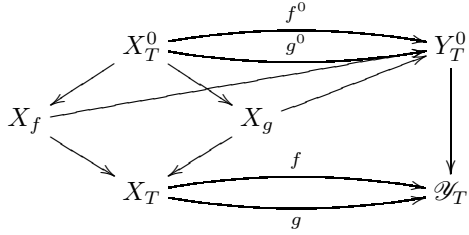
\[
 \xymatrix@R1pc{
 & X^0_T \ar[ld] \ar[rd] \ar@/^/[rrr]^{f^0} \ar@/_/[rrr]^{g^0}
  & & & Y^0_T \ar[dd]\\
 X_f \ar[rrrru] \ar[rd] & & X_g \ar[rru] \ar[ld]\\
 & X_T \ar@/^/[rrr]^{f} \ar@/_/[rrr]_g & & & \Y_T
 }
 \]
Let $X^1_T = X^0_T \times_{\X_T} X^0_T$  and  $Y^1_T = Y^0_T \times_{\Y_T} Y^0_T$.
Then $X^0_T$ and $X^1_T$ are proper and flat algebraic spaces over $T$.
Therefore the functors $\HOMs(X^0_T,Y^1_T)$, $\HOMs(X^1_T,Y^0_T)$
and $\HOMs(X^1_T,Y^1_T)$ are representable by separated algebraic spaces
over $T$. The algebraic space $\Isom_T(f,g)$ can be identified
with a closed subspace of $\HOMs(X^0_T,Y^1_T)$ whose point $\alpha$ satisfies
$p_1 \circ \alpha = f^0, p_2 \circ \alpha = g^0$ and
$\alpha \circ p_1 = \alpha \circ p_2$.
Hence $\Isom_T(f,g)$ is separated and quasicompact.
\end{proof}

\section{Application: the Picard stack}

Let $\X$ be an algebraic stack over $S$. The Picard 2-functor $\Pic_{\X}$
from the category of affine noetherian schemes over $S$
to the 2-category of groupoids is defined by
\[
 \Pic_{\X}(T) = \mbox{the category of line bundles on $\X_T$}.
\]
as in \cite[14.4.7]{LM}.
Then we have
\begin{theorem}
 If $\X$ is proper and flat over $S$, then $\Pic_{\X}$ is an algebraic
 stack in Artin's sense.
\end{theorem}
\begin{proof}
To give a line bundle on $\X$ is equivalent to
give a morphism $\X \to B\Gm$. Here $B\Gm$ denotes the classifying stack
of the multiplicative group $\Gm$.
Therefore
\[
 \Pic_{\X} = \HOMs(\X,B\Gm).
\]
This is an algebraic stack in Artin's sense by Theorem \ref{thm:main}.
\end{proof}

\appendix

\section{The semicontinuity theorem for proper algebraic stacks}

Let $x: \X \to T$ be a proper algebraic stack over
an affine scheme $T = \Spec A$
and $\F$ a coherent sheaf of $O_{\X}$-modules on $\X$.
Suppose that $T$ is reduced and $\F$ is flat over $T$. 
For each point $t$ of $T$, let $\X_t$ be the fiber over $t$ and
$\F_t = \F \tensor_{O_T} k(t)$. 

\begin{theorem}
 \label{thm:semi-con}
 $ $
\begin{enumerate}
 \item The function on $T$ defined by
       \[
       t \mapsto \dim_{k(t)}H^i(\X_t,\F_t)
	\]
       is upper semi-continuous on $Y$.
 \item There is an open subscheme $U \subset X$ in which
       \[
       R^i x_* \F \tensor_{O_T} k(t) \to H^i(\X_t,\F_t)
       \]
       is an isomorphism.
\end{enumerate}
\end{theorem} 

The proof is almost the same as one in \cite[5]{Mu}.
The key is the following lemma:
\begin{lemma}
 \label{lem:cpx_k}
 Let $\X$, $T$ and $\F$ be as above.
 For each positive integer $N$,
 there is a complex
 \[
 K^\dot : 0  \to K^0 \to K^1 \to \dots \to K^N \to 0
 \]
 of finitely generated projective $A$-modules and isomorphisms
 \[
 H^i(\X \times_T \Spec B, \F \times_A B)
 \isomto H^i(K^\dot \otimes_A B) \quad (0 < i < N)
 \]
 functorial on $A$-algebra $B$.
\end{lemma}

\begin{remark}
 This is a generalization of
 the second theorem in \cite[5]{Mu}.
 The first theorem in \cite[5]{Mu}
 which claims direct images of proper schemes are coherent
 also holds in the case of proper algebraic stacks \cite[Theorem 1]{Fa}.
 We have to limit $i<N$ because cohomological dimension
 of an algebraic stack may be infinite.
 Note that Lemma 1 and Lemma 2 in the proof of \cite[5]{Mu}
 concern only modules on $A$, and the same discussion applies
 to our case.
\end{remark}

\begin{proof}[Proof of Lemma \ref{lem:cpx_k}]
 Let $P^0 : X^0 \to \X$ be a presentation with $X^0$ affine
 and $X^\dot = \cosq{X^0}{\X}$. Then by cohomological descent,
 we have an isomorphism
 \[
  H^i(\X,\F) \simeq H^i(X^\dot, {P^\dot}^* \F).
 \]
 Since $X^0$ is affine and $\X$ is separated, $X^n$ is affine for all
 $n$ and
 $H^i(X^n,{P^n}^* \F) = 0$ for $i>0$.
 Let
 \[
  C^n = H^0(X^n,{P^n}^* \F)
 \]
 and $C^\dot$ be the alternating cochain. Then we have
 \[
  H^i(\X,\F) \simeq H^i(C^\dot).
 \]
 Note that $H^i(C^\dot)$ is a finite $A$-module because
 $\F$ is coherent.
 Moreover, for any $A$-algebra $B$,
 \[
 P^0_B : X^0_B := X^0 \times_T \Spec B \to \X \times_T \Spec B =: \X_B
 \]
 is a presentation from affine scheme and
 \[
  H^0(X^n_B,{P^n_B}^* \F \tensor_A B) \simeq
 H^0(X^n,{P^n}^*\F) \tensor_A B
 \]
 because $\F$ is flat. Therefore we have functorial isomorphisms
 \[
  H^i(\X_B,\F \tensor_A B) \simeq H^i(C^\dot \tensor_A B) \quad (i>0).
 \]

 Now replace $C^\dot$ by its truncation $\tau_{\leq N} C^\dot$
 and construct $K^\dot$ by descending induction
 as in \cite[5 Lemma 1]{Mu}. This is the desired
 complex.
\end{proof}

Fix $N$ sufficiently large. Then by Lemma \ref{lem:cpx_k},
We can reduce Theorem \ref{thm:semi-con} to statements in
homological algebra as in corollaries of \cite[5]{Mu}.
Proofs of these corollaries also works for our case.

\end{document}